\documentclass[11pt]{article}
\usepackage{ulem}
\usepackage{graphicx,amsmath,amsfonts,amssymb,color,mathrsfs, amsmath,amsthm}
\usepackage{epsfig}
\newcommand{\chuhao}{\fontsize{19pt}{\baselineskip}\selectfont}

\newcommand{\BOX}{\hfill $\Box$}

\numberwithin{equation}{section}

 \newtheorem{theorem}{Theorem}[section]
 \newtheorem{lemma}{Lemma}[section]

 \newtheorem{remark}{Remark}[section]

 \newtheorem{example}{Example}[section]

 \newtheorem{definition}{Definition}[section]

\title{\bf\color{black} \chuhao{Asymptotic behavior of the loss probability
for an M/G/1/N queue with vacations}}
\author{
Yuanyuan Liu\thanks{Corresponding author. Postal address: School
of Mathematics and Statistics, Railway Campus, Central South University,
Changsha, Hunan, 410075, China; Email address: liuyy@csu.edu.cn;
Tel.: +86 87312655267;  Fax: +86 87315586475.
\newline \indent \hspace*{-2.5mm}
$^{**}$ Postal address: School of Mathematics and Statistics,
Carleton University, 1125 Colonel By Drive, Ottawa, ON Canada K1S
5B6; E-mail address: zhao@math.carleton.ca.}, \ \ Central South
University
\\Yiqiang Q. Zhao$^{**}$,  Carleton University
}
\date{Accepted by \textit{Applied Mathematical Modelling} April-1, 2012}

\begin{document}
 \maketitle

\begin{abstract}
In this paper, asymptotic properties of the loss probability are
considered for an M/G/1/N queue with server vacations and exhaustive
service discipline, denoted by an M/G/1/N -(V, E)-queue. Exact
asymptotic rates of the loss probability are obtained for the cases
in which the traffic intensity is smaller than, equal to and greater
than one, respectively. When the vacation time is zero, the model
considered degenerates to the standard M/G/1/N queue. For this
standard queueing model, our analysis provides new or extended
asymptotic results for the loss probability. In terms of the duality
relationship between the M/G/1/N and GI/M/1/N queues, we also
provide asymptotic properties for the standard GI/M/1/N model.
\medskip


\end{abstract}

\section{Introduction}

For a queueing system with a finite system capacity $N$, the loss
probability $P_{\mbox{loss}}(N)$, defined as the steady-state
probability with which an arriving customer finds no free room in
the queue (or buffer) and is forced to leave immediately, is an
important performance measure in real applications, see e.g.
\cite{{Alfa2000},{Bai1994},{Ishizaki1999},{Miy90},{wil88}}. Simple
closed-form solutions for the loss probability are available only
for a very limited number of special cases such as the M/M/c/K+c
queue. For others, numerical computations and simulation methods are common
tools for the analysis of the loss probability, but they cannot
provide closed-form properties. Another important aspect is the
asymptotic of the loss probability, which can lead to
approximations, performance bounds among others. Also,
approximations can be often obtained through asymptotic properties
in the invariant measure for the corresponding infinite-buffer
system, for which the analysis is usually less challenging than that
for the finite-buffer one, see, e.g. \cite{{Bai1992}, {Bai1994},
 {Kim2008}, {Miy2007}}. This paper, investigating
 the asymptotics of the loss probability for an M/G/1/N queue with
 vacations, will go along this direction.

Now we introduce this system based on the standard M/G/1/N queue.
Recall that in the M/G/1/N queue, the arrival process is a Poisson
process with rate $\lambda$ and the service times are i.i.d. random
variables with a generic random variable $S$, which is independent
of the arrival process, and has distribution function $S(x)$ and
probability density function $s(x)$. The queue has a finite buffer
of size $N-1$ to store the incoming customs. Therefore, $N$ is the
maximal system capacity for holding customers including the one in
service. To introduce the vacation model M/G/1/N-(V, E)-queue, we
further assume in the M/G/1/N queue that whenever the system becomes
empty, the server starts a vacation of random length $V$
immediately. If, when the vacation ends, the system is still empty,
the server takes another independent vacation of the same length as
in the previous one until there is at least one customer waiting in
the buffer upon return of the server from a vacation. From then, the
server keeps serving customers until the system becomes empty (all
customers have been served or exhausted). The vacation time $V$ is
assumed to be a non-negative random variable with the distribution
function $V(x)$.  Then, $\rho=\lambda E[S]$ is the traffic intensity
of the model. Throughout the paper, assume that $E[S]<\infty$ and
$E[V]<\infty$. Let $\mathbb{Z}_+$ and $\mathbb{C}$ be the sets of
the non-negative integers and the complex numbers, respectively.
Denote by $a_j$ and $\nu_j$ the probability that $j$ customers join
the queue during a service time $S$ and a vacation time $V$,
respectively, i.e., for $j\in \mathbb{Z}_+$,
\[
    a_j=\int_0^\infty \frac{(\lambda t)^j}{j!}e^{-\lambda t}dS(t),
\]
and
\[ \nu_j=\int_0^\infty \frac{(\lambda t)^j}{j!}e^{-\lambda t}dV(t).
\]

 Various performance measures, such as  the loss probability, the queue length, the waiting time, the busy period, the system throughput and so on, were investigated for the M/G/1/N-(V, E)-queue  in the literature,  e.g.  \cite{{Lee84},{keilson1989},{Takagi94},{Fr97}}.
The loss probability, which was studied in \cite{{Lee84},{keilson1989},{Fr97}},  has been treated as a central issue among these performance measures. In these works, the loss probability for the finite M/G/1/N-(V, E)-queue was expressed in terms of
the invariant probabilities for the embedded queue of the finite or infinite buffer queues, which are again
unknown. As far as we know, for the model studied here no asymptotic analysis,
as $N \to \infty$, for the loss probability  has been reported, which stimulated us to  provide a detailed characterization of the asymptotic analysis of the loss probability.  In
the literature, asymptotic results for the M/G/1/N queue (without
server vacations) can be found, for example, in  \cite{Bai1992}. In
addition, asymptotic rates for the loss probability of the GI/M/1/N
queue were obtained in (\cite{{Bai1992},{Kim2000}}), by using the
duality relationship between the M/G/1/N queue and the GI/M/1/N
queue. However, asymptotic analysis of the loss probability seems not complete even for the M/G/1/N and GI/M/1/N
queues and some important cases remain unexplored, which stimulated us to introduce  new asymptotic methods to study the loss probability for the more general  M/G/1/N-(V, E)-queue  so that the methods or the results can be applied to derive new asymptotic behavior for those queues.

In this paper we will provide a detailed characterization of asymptotic behavior of the loss probability for the
M/G/1/N-(V,E)-queue. Some of the results are intuitively expected
and others are not, and some of results are easy consequences from
the existing literature and others deserve more detailed analysis.
The characterization is obtained by investigating the dominant
singular point in terms of singularity analysis, such as the
Tauberian theorem, Tauberian-like theorem (Theorem A.1 in Appendix) or a direct complex analysis. As special cases,
for the standard M/G/1/N queue and the GI/M/1/N queue, we will
provide new asymptotic results. Specifically, our main contributions
include: (1) For $\rho <1$ and under a light-tailed condition, two
cases are studied in Theorem~\ref{the:2.2}: for the first case, the
decay is light but not exactly geometric, a new asymptotic behavior,
and for the second case, decay is exactly geometric; (2) Also for
$\rho <1$, asymptotic results when either the service or the
vacation variable is heavy-tailed are provided in
Theorem~\ref{the:2.1}. This case was not considered in
\cite{Bai1992};
 (3) For $\rho=1$,
asymptotics are studied in Theorem~\ref{the:2.3} under the
assumption that the second moment of the service time is finite and
infinite, respectively. The infinite case was not considered in
\cite{Bai1992}; (4) An asymptotic result for $\rho>1$ is reported in
Theorem~\ref{the:2.4}, which can be expected from the corresponding
result in \cite{Bai1992}; (5) When $V=0$, the vacation model is
degenerated to the standard M/G/1/N queue, for which asymptotic
results are summarized in Theorem~\ref{the:3.1}. We also provide
asymptotic results for the GI/M/1/N queue in Theorem~\ref{the:3.2}
by using the dual relationship with the M/G/1/N queue. These two
theorems provide extensions of existing literature results.

The rest of the paper is organized into: Section~2, which bridges
the asymptotics of the loss probability of the finite queue and
decay behavior of the invariant measure of the corresponding
infinite queue; Section~3, which is devoted to asymptotic analysis
of the M/G/1/N-(V,E)-queue and contains main results of the paper;
Section~4, in which asymptotic results for the standard M/G/1/N
queue the GI/M/1/N queue are provided; Section~5, which provides concluding remarks; and an Appendix, which
contains some properties in asymptotic analysis used in this paper
for convenience.

\section{Connection with the  M/G/1/$\infty$-(V, E)-queue}

Let $L(t)$ be the queue length of the  M/G/1/N-(V,E)-queue at time
$t$. Since the state space is finite, the system $L(t)$ is always
stable. It is known that $L(t)$ itself is not a Markov process, but
it can be analyzed through its embedded chain $L(\tau_n)$ at the
departure epochs $\tau_n$, $n\in \mathbb{Z_+}$. The transition
matrix of the embedded chain $L(\tau_n)$ is given by
\[
 P(N)=\left (
 \begin{array}{cccccc}
   b_0 & b_1 & b_2  & ...& b_{N-2} & 1-\sum_{k=0}^{N-2}b_k  \\
   a_0 & a_1 & a_2  & ... &a_{N-2} & 1-\sum_{k=0}^{N-2}a_k\\
   0   & a_0 & a_1  & ... & a_{N-3} & 1-\sum_{k=0}^{N-3}a_k \\
   0   & 0   & a_0  & ...   & a_{N-4}    & 1-\sum_{k=0}^{N-4}a_k \\
   \vdots & \vdots & \vdots &\ddots &\vdots &\vdots \\
   0 & 0 & 0 & ...& a_0    &    1-a_0
   \end{array}\right )_{N\times N},
\]
where $b_j=\sum_{i=1}^{j+1}\frac{\nu_i}{1-\nu_0}a_{j-i+1}$, $j\in
\mathbb{Z_+}$. Obviously, $P(N)$ is positive recurrent. Let $\pi(N)$
be the probability invariant vector of $P(N)$, i.e. $\pi(N)
P(N)=\pi(N)$, or
\begin{equation}\label{2.1}
  \pi_j(N)=\pi_0(N) b_j+\sum_{k=1}^{j+1}\pi_k(N) a_{j+1-k},\ \ j=0,\ldots, N-2,
\end{equation}
and $ \sum_{k=0}^{N-1}\pi_k(N)=1$.

Let $\pi_j^*(N)=\lim_{t\rightarrow \infty} P[L(t)=j]$, $j=0,1, \ldots,
N$. By the Poisson Arrival See Time Average Property, we see that
$\pi_j^*(N)$ is also the probability that there are $j$ customers in
the system just before an arrival, which implies that
$P_{\mbox{loss}}(N)=\pi_N^*(N)$. By Theorem~4.2 in \cite{Fr97}, we have
\begin{equation}\label{2.2}
  P_{\mbox{loss}}(N)=1-\frac{(1-\nu_0)
  \lambda^{-1}}{E[V]\pi_0(N)+E[S](1-\nu_0)}.
\end{equation}
Therefore, the asymptotic analysis of the loss probability solely
depends on the property of $\pi_0(N)$. In the following we establish
a relationship between the stationary probabilities for the
finite-buffer queue and the invariant measure  for the corresponding
infinite-buffer queue for all values of the traffic intensity
$\rho$.

To connect the loss probability $P_{\mbox{loss}}(N)$ with the
invariant measure of the corresponding infinite queue, we now
introduce some related notations. For sequences $\{x_n, n\in
\mathbb{Z}_+\}$ and $\{y_n, n\in \mathbb{Z}_+\}$ of non-negative
real numbers, define
\[
    S_{x}(n)=\sum_{j=0}^{n-1} x_j, \ \ \bar{S}_{x}(n)=\sum_{j=n}^\infty x_j,\ \ n\geq 1,
\]
and write $x_n\sim y_n$ to imply $\lim_{n\rightarrow \infty}
{x_n}/{y_n}=1$. Similarly, for two functions  $f(x)$ and $g(x)$,
$f(x) \sim g(x)$ as $x\rightarrow x_0-$ stands for
$\lim_{x\rightarrow x_0- } f(x)/g(x)=1$, and $f(x) \sim g(x)$ for
$\lim_{x\rightarrow +\infty } f(x)/g(x)=1$. For any distribution
function $F(x)$ on the nonnegative real line, let $f(x)$ be its
probability density function when it exists,
$\overline{F}(x)=1-F(x)$ be its complementary distribution function,
\[
F^*(t)=\int_0^\infty e^{-t x}dF(x)=\int_0^\infty e^{-t x}f(x) dx,\ \
t\in \mathbb{C}
\]
be its Laplace transform, and
\[
  F_e(x)=\frac{1}{E[X]} \int_0^x \overline{F}(u) du
\]
be its  equilibrium distribution whenever $E[X]$ is finite. The
distribution function $F(x)$ is said to be regularly varying if
$\overline{F}(x)=x^{-\alpha} L(x)$, where $L(x)$ is a slowly varying
function that  satisfies $\lim_{x\rightarrow + \infty} \frac{L(t
x)}{L(x)}=1$ for any $t>0$.

The corresponding infinite-buffer vacation model is referred to as
the M/G/1/$\infty$-(V, E)-queue, and its embedded chain has the
following transition matrix:
\[
    P=\left (
 \begin{array}{ccccc}
   b_0 & b_1 & b_2 & b_3 & \cdots \\
   a_0 & a_1 & a_2 & a_3 &\cdots \\
   0  & a_0 & a_1 & a_2 & \cdots \\
   0  & 0 & a_0 & a_1 & \cdots \\
   \cdots &\cdots& \cdots & \cdots& \cdots \
   \end{array}\right ),
   \]
where the coefficients $b_i$ and $a_i$ are the same as those in
$P(N)$. Regardless of stability  or non-stability  of the
infinite-buffer queue, $P$ has a unique invariant measure, up to a
multiplicative constant, due to its special structure.  Let $\pi$ be
the invariant measure of $P$, i.e. $\pi P=\pi$, or
\begin{equation}\label{2.3}
\pi_j =\pi_0 b_j + \sum_{k=1}^{j+1}\pi_k a_{j+1-k}, \ \ j\geq 0.
\end{equation}
Define the generating functions
 \[
    B(z)=\sum_{k=0}^\infty b_k z^k,\ \ A(z)=\sum_{k=0}^\infty a_k
 z^k\ \ \mbox{and} \ \ \Pi(z)=\sum_{k=0}^\infty \pi_k z^k
 \]
 for $b_k$, $a_k$ and  $\pi_k$,
 respectively.
 It is routine to calculate
 \[
   A(z)= S^*(\lambda-\lambda z),\ \ B(z)=\frac{S^*(\lambda-\lambda z)[V^*(\lambda-\lambda
   z)-\nu_0]}{z(1-\nu_0)}.
 \]\
 From (\ref{2.3}), we have
\[
 \Pi(z)[z-A(z)]=\pi_0 [z B(z)- A(z)],
\]
whenever $\Pi(z)$ is finite, that is,
\begin{equation}\label{2.35}
 \Pi(z)[z-S^*(\lambda-\lambda z)]=\frac{\pi_0[V^*(\lambda-\lambda z)-1] S^*(\lambda-\lambda z)}{1-\nu_0}.
\end{equation}

Observing that all equations of (\ref{2.1}) are identical to the
corresponding ones in (\ref{2.3}), we have
\begin{equation}\label{2.4}
\pi_i(N)= \frac{\pi_i}{S_{\pi}(N)}, \ \ i=0,1, \ldots,
    N-1.
\end{equation}

\begin{remark}
Note that $P$ is positive recurrent, or $\pi$ is a probability
measure, if and only if $\rho<1$.
\end{remark}

The following lemma  will later be frequently used to analyze the
asymptotic behavior of the loss probability in Section 3.
\begin{lemma} \label{lem:2.1}
For the case of $\rho<1$, we have
\begin{equation}\label{2.5}
 \lim_{N\rightarrow \infty} \frac{1}{{\overline S}_{\pi}(N)}P_{\mbox{loss}}(N) =1-\rho;
\end{equation}
and for the case of $\rho\geq 1$, we have
\begin{equation}\label{2.6}
 \lim_{N\rightarrow \infty}S_{\pi}(N)\left[ P_{\mbox{loss}}(N) - \left (1-\frac{1}{\rho} \right )\right]= \frac{ E[V] \pi_0} {\rho
 E[S](1-\nu_0)}.
\end{equation}
\end{lemma}

\proof  Inserting (2.4) into (2.2) yields
\begin{equation}\label{2.7}
  P_{\mbox{loss}}(N)=\frac{E[V]\pi_0 +
(E[S]-\lambda^{-1})(1-\nu_0)S_{\pi}(N)} {E[V]\pi_0 + E[S](1-\nu_0)
S_{\pi}(N)}.
\end{equation}

When $\rho<1$, we have $\sum_{i=0}^\infty \pi_i=1$ and $
\pi_0=\frac{(1-\rho)(1-v_0)}{\lambda E[V]}$. Use (\ref{2.7}) to
derive
\[
 P_{\mbox{loss}}(N) =\frac{(1-\rho) {\overline S}_{\pi}(N)}{(1-\rho)+\rho
S_{\pi}(N)},
\]
from which (\ref{2.5}) follows easily.

When $\rho\geq 1$, we have from (\ref{2.7})
\[
 \lim_{N\rightarrow \infty}
 P_{\mbox{loss}}(N)=1-\frac{1}{\rho},
\]
since $S_{\pi}(N) \to \infty$ as $N \to \infty$. Write (\ref{2.7})
as
\[
 P_{\mbox{loss}}(N) - \left (1-\frac{1}{\rho} \right )= \frac{ E[V] \pi_0} {\rho
 E[V]\pi_0+\rho S_{\pi}(N) E[S](1-\nu_0)}.
\]
Then, (\ref{2.6}) follows immediately. \BOX

 From Lemma~\ref{lem:2.1}, we
see that the asymptotic behavior of the loss probability
$P_{\mbox{loss}}(N)$ as $N\rightarrow \infty$ relies on the decay
property in the invariant measure $\pi$. For a stable system, we
study the decay rate of the tail probability ${\overline
S}_{\pi}(N)$, and for a non-stable system, we analyze the decay of
$1/S_{\pi}(N)$.

For the Laplace transform $F^*(\lambda-\lambda z)$
($F^*(\lambda-\lambda z)$ $=F^*(t)|_{t=\lambda-\lambda z}$), define
$R_{F^*} \geq 1$ to be the leftmost singular point of the functions
$F^*(\lambda-\lambda z)$,
\[
  {F^*}'(\lambda-\lambda z)= \frac{d F^*(\lambda-\lambda
u)}{d u}|_{u=z}=\int_0^\infty \lambda t e^{- (\lambda-\lambda
z)t}f(t)dt,
\]
and
\[
  {F^*}''(\lambda-\lambda z)= \frac{d^2 F^*(\lambda-\lambda
u)}{d u^2}|_{u=z}=\int_0^\infty (\lambda t)^2 e^{- (\lambda-\lambda
z)t}f(t)dt.
\]

Observe that ${S^*}(\lambda-\lambda\cdot 1)=1$,
${S^*}'(\lambda-\lambda\cdot 1)=\rho$ and $S^*(\lambda-\lambda z)$
is a strictly increasing and convex function about $z$ when $z\geq
0$. This observation yields the following lemma, which will later be
used  to investigate the decay of ${\overline S}_{\pi}(N)$ and
$1/S_{\pi}(N)$.

\begin{lemma}
(i) If $\rho<1$ and $1<R_{S^*}< S^*(\lambda-\lambda R_{S^*})$, then
there exists at most one solution $z_1$ to the equation
$z=S^*(\lambda-\lambda z)$, $1<z<R_{S^*}$.

(ii) If $\rho>1$, then there exists exactly one solution $z_2$ to
the equation $z=S^*(\lambda-\lambda z)$, $0<z<1$.

\end{lemma}

\section{ Asymptotic behavior for the  M/G/1/N-(V, E)-queue}

We perform the asymptotic analysis according to three cases: (i)
$\rho<1$; (ii) $\rho=1$;  and  (iii) $\rho>1$.


\subsection{The case when $\rho<1$}

For this case, we investigate the asymptotic behavior for two
subcases according to $\min\{R_{S^*}, R_{V^*}\}$\\$>1$ or
$\min\{R_{S^*}, R_{V^*}\}=1$. The former subcase corresponds to
light-tailed characterizations, for which the asymptotic rates are
either exactly geometric or not exactly geometric. The latter
subcase corresponds to that one of the sequences $\{b_i\}$ and
$\{a_i\}$ is heavy-tailed, for which a regular varying condition is
imposed on the service time and vacation time distributions.

\begin{theorem} \label{the:2.2}  Suppose that  $\rho<1$ and  $r:=\min\{R_{S^*},
R_{V^*}\}>1$.
\begin{description}
\item[(i)]  Suppose that $r=R_{V^*}$ and $S^*(\lambda-\lambda r)<r$.
If $V^*(\lambda-\lambda z)$ can be analytically extended to a
$\Delta(r)$-domain  (see Definition A.1 in Appendix) and for some
$\theta>0$, $
 V^*(\lambda-\lambda z) \sim \frac{c}{(r-z)^\theta}, \ \mbox{as} \ \ z\rightarrow r, z\in \Delta(r)
$, then
   \begin{equation}\label{2.12}
  \lim_{N\rightarrow \infty} N^{1-\theta}r^{N} P_{\mbox{loss}}(N)
  =\frac{c (1-\rho)^2
S^*(\lambda-\lambda r)}{ r^{\theta-1} \Gamma(\theta)\lambda E[V]
(r-1) [r-S^*(\lambda-\lambda r)]}.
\end{equation}

\item[(ii)]
 Suppose that $r=R_{s^*}$ and $r< S^*(\lambda-\lambda
r)$, then there exists exactly one  solution $z_1$ to the equation
$z=S^*(\lambda-\lambda r)$,  $1<z<r$, such that
\begin{equation}\label{2.11}
  \lim_{N\rightarrow \infty} z_1^{N} P_{\mbox{loss}}(N)
  =\frac{ z_1(1-\rho)^2
  [1-V^*(\lambda-\lambda z_1)]}{\lambda E[V](z_1-1)[1- {S^*}'(\lambda-\lambda
  z_1)]}.
\end{equation}

\end{description}
 \end{theorem}
\proof (i) Due to the assumption, we know that  $V^*(\lambda-\lambda
z)$ is analytic in a $\Delta(r)$-domain, which implies that there
exists some $r_1$, $r_1>r$ and some $\alpha$,
$0<\alpha<\frac{\pi}{2}$ such that $V^*(\lambda-\lambda z)$ is
analytic in the domain $\Delta(r, \alpha, r_1)$. Thus we can find an
$r_2$ such that $r<r_2<r_1$ and $V^*(\lambda-\lambda z)$ is analytic
in the domain $\Delta(r, \alpha, r_2)$. Hence by the expression
(\ref{2.35}), we know that $\Pi(z)$ can be analytically extended to
this domain $\Delta(r, \alpha, r_2)$, i.e. $\Pi(z)$ is
$\Delta(r)$-analytic.

 Since $V^*(\lambda-\lambda z) \sim
\frac{c}{(r-z)^\theta}$ as $ z \rightarrow r$ in the complex domain
$\Delta(r, \alpha, r_2)$, we have
\begin{eqnarray*}
\lim_{z\rightarrow r} \Pi(z)(r-z)^\theta = \frac{c (1-\rho)
S^*(\lambda-\lambda r)}{\lambda E[V][r-S^*(\lambda-\lambda r)]},
\end{eqnarray*}
that is,
\[
  \sum_{n=0}^\infty \pi_n r^n (\frac{z}{r})^n \sim \frac{c (1-\rho)
S^*(\lambda-\lambda r)}{\lambda E[V][r-S^*(\lambda-\lambda r)]}
\frac{1}{r^\theta (1-\frac{z}{r})^\theta}, \ \ \mbox{as}\ \ z
\rightarrow r, z\in \Delta(r, \alpha, r_2).
\]
It follows essentially from a result in \cite{Fla09} (see
Theorem~\ref{the:4.1} in Appendix) that
\[
  \pi_n \sim \frac{c (1-\rho)
S^*(\lambda-\lambda r)}{ r^\theta \Gamma(\theta)\lambda  E[V]
[r-S^*(\lambda-\lambda r)]} n^ {\theta-1} r^{-n},
 \]
from which and the Stolz-Ces\'{a}ro theorem, we have
\[  \overline{S}_\pi(n) \sim   \frac{c (1-\rho) S^*(\lambda-\lambda r)}{ r^{\theta-1} \Gamma(\theta)\lambda E[V]
(r-1) [r-S^*(\lambda-\lambda r)]} n^ {\theta-1}   r^{-n}. \] Then,
assertion (\ref{2.12}) follows from Lemma~\ref{lem:2.1} immediately.

(ii) It follows from Lemma 2.2 that there exists exactly one
solution $z_1$ to the equation $z=S^*(\lambda-\lambda z)$, $1<z<r$.
Since
\[
\lim_{z\rightarrow z_1-} \Pi(z)(z_1-z)
=  \frac{z_1 (1-\rho)[1-V^*(\lambda-\lambda z_1) ]}{\lambda
E[V][1-{S^*}'(\lambda-\lambda z_1)]},
\]
 $z_1$ is a simple pole for $\Pi(z)$. For any $z\in \mathbb{C}$ such that $|z|=z_1$ and $z\neq z_1$, we have
\[
  |S^*(\lambda-\lambda z)|\leq  \int_0^\infty e^{-\lambda x} |e^{\lambda z x}| d
  S(x)
  < S^*(\lambda-\lambda |z|)=z_1,
\]
which implies that there is only one solution $z=z_1$ to the
equation $z=S^*(\lambda-\lambda z)$ on the circle $|z|=z_1$. Hence
$z=z_1$ is the only singular point of $\Pi(z)$ on the circle
$|z|=z_1$.  By Theorem~5.2.1 in \cite{wilf94}, a standard result on
asymptotics of complex functions, we have
\[
\pi_n\sim \frac{(1-\rho)[1-V^*(\lambda-\lambda z_1) ]}{\lambda
E[V][1-{S^*}'(\lambda-\lambda z_1)]}z_1^{-n},
\]
from which and the  Stolz-Ces\'{a}ro theorem, it follows that
\begin{equation*}
\overline{S}_{\pi}(n)\sim \frac{z_1(1-\rho)[1-V^*(\lambda-\lambda
z_1) ]}{\lambda E[V](z_1-1)[1-{S^*}'(\lambda-\lambda z_1)]}z_1^{-n}.
\end{equation*}
The assertion (\ref{2.11}) follows from Lemma~\ref{lem:2.1}
immediately. \BOX

\begin{remark} (i) In the above theorem,  the assumption
$S^*(\lambda-\lambda r)<r$ in (i), which is equivalent to that there
is no solution to the equation $z=S^*(\lambda-\lambda z)$ for
$1<z<r$, ensures that the dominant singular point of $\Pi(z)$ is
determined by $V^*(\lambda-\lambda z)$, while the assumption in (ii)
ensures that the dominant singular point of $\Pi(z)$ is determined
by the solution to the equation $z=S^*(\lambda-\lambda z)$ for
$1<z<r$. (ii) In the first case, the exact asymptotic behavior of
$P_{\mbox{loss}}(N)$ reveals a non-geometric phenomenon when
$\theta\neq 1$, which is faster than the geometric rate $r^{-N}$
when $\theta<1$ and slower than the geometric rate $r^{-N}$ when
$\theta>1$. This is a different characterization from the well-known
exact geometric decay in the related literature.
\end{remark}

\begin{example}
We now give an example to illustrate that how to verify the analytic
condition in (i) of Theorem 3.1. Assume that (i): $\lambda=1$, $V$
has probability density function $v(x)=\frac{\alpha^{k+1} x^k}{k!}
e^{-\alpha x}, \alpha>0,x\geq 0$, and $S$ has  probability density
function $s(x)=pe^{-px}, p>0, x\geq 0$. Then it is easy to calculate
\[
 E[V]=\frac{k+1}{\alpha}, \ \ E[S]=\frac{1}{p}, \ \ \rho=\lambda
 E[S]=1/p,
\]
\[
  V^*(\lambda-\lambda z)=V^*(1-z)=\int_0^{\infty}
  e^{-(1-z)x}v(x)dx=\frac{\alpha^{k+1}}{(1+\alpha-z)^{k+1}},\ \
  Re(z)<1+\alpha,
\]
and
\[
   S^*(\lambda-\lambda z)=S^*(1-z)=\int_0^{\infty}
  e^{-(1-z)x}s(x)dx=\frac{p}{1+p-z},\ \ Re(z)<1+p.
 \]
So, we have $R_{V^*}=1+\alpha$, $R_{S^*}=1+p$. Assume that (ii):
$p>1+\alpha$, then we have $r=R_{V^*}$ and $S^*(\lambda-\lambda
r)<r$. Obviously, the function
$\frac{\alpha^{k+1}}{(1+\alpha-z)^{k+1}}$ is analytic in
$\mathbb{C}\setminus \{1+\alpha\}$, i.e. the whole complex plane
except for the point $z=1+\alpha$. Hence, the Laplace transform
$V^*(\lambda-\lambda z)$ can be analytically extended to
$\mathbb{C}\setminus \{1+\alpha\}$, and the $\Delta(r)$-analytic
condition holds automatically. Thus, under the assumptions (i) and
(ii), every condition in (i) of Theorem 3.1 is satisfied, so we have
\[
  \lim_{N\rightarrow \infty} N^{-k} (1+\alpha)^{N} P_{\mbox{loss}}(N)
  =\frac{ \alpha^{k}p(1-1/p)^2}{  (1+\alpha)^k (k+1)!  (p-1-\alpha)}.
\]
\end{example} \BOX

\begin{theorem} \label{the:2.1} Assume that  $\rho<1$ and $\min\{R_{S^*},
R_{V^*}\}=1$.
 \begin{description}
\item[(i)] If $\overline{V}(x)\sim c \overline{S}(x)\sim x^{-\alpha} L(x)$
for some $c>0$ and $\alpha>1$, then
\begin{equation}\label{2.8}
  \lim_{N\rightarrow \infty} \frac{N^{\alpha-1}}{L(\frac{N}{\lambda})} P_{\mbox{loss}}(N)
  =\left(\frac{1-\rho}{E[V]}+\frac{\rho}{cE[S]}\right)\frac{\lambda^{\alpha-1}}{\alpha-1}.
\end{equation}
\item[(ii)] If $\overline{V}(x)=o(\overline{S}(x))$ and $\overline{S}(x)\sim x^{-\alpha}
L(x)$ for some $\alpha>1$, then
\begin{equation}\label{2.9}
  \lim_{N\rightarrow \infty} \frac{N^{\alpha-1}}{L(\frac{N}{\lambda})} P_{\mbox{loss}}(N)
  =\frac{\lambda^\alpha}{\alpha-1}.
\end{equation}
\item[(iii)] If $\overline{S}(x)=o(\overline{V}(x))$ and $\overline{V}(x)\sim x^{-\alpha}
L(x)$ for some $\alpha>1$, then
\begin{equation}\label{2.10}
  \lim_{N\rightarrow \infty} \frac{N^{\alpha-1}}{L(\frac{N}{\lambda})} P_{\mbox{loss}}(N)
  =\frac{(1-\rho) \lambda^{\alpha-1}}{(\alpha-1)E[V]}.
\end{equation}
\end{description}
\end{theorem}

\proof (i)
 It was shown in Proposition~1.5.10 in \cite{Bing1987} that if $\overline{F}(x)\sim
x^{-\alpha} L(x)$ for some $\alpha>1$, then the  equilibrium
distribution
\begin{equation}\label{2.103}
 \overline{F}_e(x)\sim
\frac{x^{-(\alpha-1)}L(x)}{E[X](\alpha-1)}.
\end{equation}
Using L'Hospital rule obtains
\begin{equation}\label{2.105}
 \lim_{x\rightarrow \infty} \frac{\overline{V}_e(x)}{
  \overline{S}_e(x)}= \lim_{x\rightarrow \infty} \frac{E[S](E[V]-\int_0^x \overline{V}(t)dt)}{E[V](E[S]-\int_0^x \overline{S}(t)dt)}
  =  \frac{E[S]}{E[V]}\lim_{x\rightarrow \infty}\frac{\overline{V}(x)}{
  \overline{S}(x)},
\end{equation}
 whenever the limit $\lim_{x\rightarrow \infty}\frac{\overline{V}(x)}{
  \overline{S}(x)}$ exists. Using (\ref{2.103}) and (\ref{2.105}) yields
\[
  \overline{V}_e(x)\sim \frac{c E[S]}{E[V]}\overline{S}_e(x)\sim
  \frac{x^{-(\alpha-1)}L(x)}{E[V](\alpha-1)}.
\]
Then by Proposition 6.4 in \cite{Asmu99}, we have
\[
  \overline{S}_{\pi}(k)\sim
  \left(\frac{1}{E[V]}+\frac{\rho}{c(1-\rho)E[S]}\right)\frac{\lambda^{\alpha-1}}{\alpha-1}k^{-(\alpha-1)}
  L\left (\frac{k}{\lambda}\right ),
\]
from which and Lemma~\ref{lem:2.1}, (\ref{2.8}) follows.

(ii) It follows from (\ref{2.103}) and (\ref{2.105}) that
\[
 \overline{S}_e(x)\sim
\frac{x^{-(\alpha-1)}L(x)}{E[S](\alpha-1)}, \ \ \overline{V}_e(x)=o(
  \overline{S}_e(x)).
\]
 From Proposition 6.1 in \cite{Asmu99}, we derive
  \[
  \overline{S}_{\pi}(k)\sim
  \frac{\lambda^\alpha}{(1-\rho)(\alpha-1)} k^{-(\alpha-1)}L\left (\frac{k}{\lambda}\right ),
 \]
from which and Lemma~\ref{lem:2.1}, (\ref{2.9}) follows.

(iii)  It follows from (\ref{2.103}) and (\ref{2.105}) that
\[
 \overline{V}_e(x)\sim
\frac{x^{-(\alpha-1)}L(x)}{E[V](\alpha-1)}, \ \ \overline{S}_e(x)=o(
  \overline{V}_e(x)).
\]
 From Proposition 6.2 in \cite{Asmu99}, we derive
  \[
  \overline{S}_{\pi}(k)\sim
  \frac{\lambda^{\alpha-1}}{(\alpha-1)E[V]}k^{-(\alpha-1)}L\left (\frac{k}{\lambda}\right ),
 \]
from which and Lemma~\ref{lem:2.1}, (\ref{2.10}) follows.
 \BOX

\begin{remark} (i) Under the assumption that $ \overline{S}(x)\sim x^{-\alpha} L(x)$
($\overline{V}(x)\sim x^{-\alpha} L(x)$), it is natural to restrict
that $\alpha>1$ due to $E[S]<\infty$ ($E[V]<\infty$), the basic
assumptions of this model. (ii) Heavy-tailed asymptotic results were
obtained for M/G/1-type Markov chains (e.g. \cite{Li2005b,
Takine2004}), in terms of asymptotics of $a(n)$ and $b(n)$, which
can not trivially lead to asymptotic results in terms of queueing
primitives or service/arrival parameters, for example, results
presented in this paper.
\end{remark}
\medskip

\subsection{The case when $\rho=1$}

For this case, if the second moment of the service time is finite,
we provide the asymptotic decay for the loss probability. When the
second moment is infinite, we consider a family of service
distributions to understand the asymptotic behavior.

\begin{theorem} \label{the:2.3} Assume that  $\rho=1$.

\begin{description}

\item[(i)] If $E[S^2]<\infty$, then
\begin{equation}\label{2.13}
  \lim_{N\rightarrow \infty} N P_{\mbox{loss}}(N)=\frac{\lambda^2 E[S^2]}{2}.
\end{equation}

\item[(ii)] If $E[S^2]=\infty$, let $s(x)\sim c x^{-\theta}$ for
$2<\theta\leq 3$, where $c$ is a constant. Then, for $2<\theta<3$,
\begin{equation}\label{2.14}
 \lim_{N\rightarrow \infty} N^{\theta-2} P_{\mbox{loss}}(N)=\frac{c \lambda^{\theta-1}
  \Gamma(\theta-1)\Gamma(4-\theta)}{(1-\theta)(2-\theta)(3-\theta)},
\end{equation}
and for $\theta =3$,
\begin{equation}\label{2.141}
 \lim_{N\rightarrow \infty} \frac{N}{\ln N} P_{\mbox{loss}}(N)=\frac{c
 \lambda^2 }{2}.
\end{equation}
\end{description}
\end{theorem}

\proof (i) Since
\[
  \lim_{z\rightarrow 1-}\frac{z-S^*(\lambda-\lambda z)}{(1-z)^2}=\lim_{z\rightarrow
  1-}\frac{1-{S^*}'(\lambda-\lambda z)}{-2(1-z)}=\frac{-{S^*}''(\lambda-\lambda \cdot 1)}{2}=-\frac{\lambda^2 E[S^2]}{2},
\]
 we have
\begin{equation*}
 \lim_{z\rightarrow 1-}\Pi(z) (1-z)= \frac{- 2\pi_0}{(1-\nu_0)\lambda^2 E[S^2]}\lim_{z\rightarrow
 1-} \frac{V^*(\lambda-\lambda z)-1}{1-z}
= \frac{2\pi_0 \lambda E[V]}{(1-\nu_0)\lambda^2 E[S^2]}.
\end{equation*}
 It follows from a Tauberian theorem (Theorem~\ref{the:4.3} in the Appendix) that
\[
  S_{\pi}(n) \sim  \frac{2\pi_0 \lambda E[V] n}{(1-\nu_0)\lambda^2
  E[S^2]}.
\]
From Lemma~\ref{lem:2.1}, we have
\[
  \lim_{N\rightarrow \infty} N P_{\mbox{loss}}(N)=\frac{\lambda E[S^2]}{2 \rho
  E[S]},
\]
which is (\ref{2.13}) by noting that $\rho=\lambda E[S]=1$.
\medskip

(ii) Let $U(x)=\int_0^x (\lambda t)^2 s(t)dt$. There are two
cases.

(1) If $2< \theta< 3$, then we have
\[
  \lim_{x\rightarrow \infty} \frac{U(x)}{x^{3-\theta}}=\frac{c \lambda^2}{3-\theta}.
\]
Obviously,
\[
 { S^*}''(\lambda-\lambda z)=\int_0^\infty (\lambda t)^2 e^{- (\lambda-\lambda
z)t}s(t)dt=\int_0^\infty e^{-
  (\lambda-\lambda z)t}dU(t).
\]
Using a Tauberian theorem (Theorem~\ref{the:4.2} in the Appendix)
yields
\[
 { S^*}''(\lambda-\lambda z) \sim \frac{c
 \lambda^{\theta-1}\Gamma(4-\theta)}{3-\theta}(1-z)^{-(3-\theta)},\ \ z\rightarrow 1-.
 \]
By L'Hospital rule, we have
\[
  \lim_{z\rightarrow 1-}\frac{z-S^*(\lambda-\lambda z)}{(1-z)^{-(1-\theta)}}=\lim_{z\rightarrow
  1-}\frac{1-{S^*}'(\lambda-\lambda z)}{-(1-\theta)(1-z)^{-(2-\theta)}}=\frac{- c
  \lambda^{\theta-1}\Gamma(4-\theta)}{(1-\theta)(2-\theta)(3-\theta)},
\]
from which it follows that
\begin{eqnarray*}
 \lim_{z\rightarrow 1-}\Pi(z) (1-z)^{-(2-\theta)}
&=& \frac{\pi_0}{(1-\nu_0)} \lim_{z\rightarrow
 1-} \frac{[V^*(\lambda-\lambda z)-1](1-z)^{-(2-\theta)}}{{\frac{-c
  \lambda^{\theta-1}\Gamma(4-\theta)}{(1-\theta)(2-\theta)(3-\theta)} (1-z)^{-(1-\theta)}}}\nonumber\\
  &=&  \frac{ \pi_0(1-\theta)(2-\theta)(3-\theta)\lambda E[V]}{c(1-\nu_0)
  \lambda^{\theta-1}\Gamma(4-\theta)}.
\end{eqnarray*}
Now, by using another Tauberian theorem (Theorem~\ref{the:4.3}),
we have
\[
  S_{\pi}(n) \sim \frac{\pi_0(1-\theta)(2-\theta)(3-\theta)\lambda E[V]}{c(1-\nu_0)
  \lambda^{\theta-1}\Gamma(4-\theta)}
  \frac{n^{\theta-2}}{\Gamma(\theta-1)},
\]
from which and Lemma~\ref{lem:2.1}, (\ref{2.14}) follows
immediately.
\medskip

(2) If $\theta =3$, then
\[
  \lim_{x\rightarrow \infty} \frac{U(x)}{\ln x}= c \lambda^2.
\]
Using a Tauberian theorem (Theorem~\ref{the:4.2})  yields
\[
 { S^*}''(\lambda-\lambda z) \sim  c
 \lambda^2 \ln \left (\frac{1}{\lambda-\lambda z}\right )=-c \lambda^2 \ln(\lambda-\lambda z).
 \]
Then, by L'Hospital rule, we obtain
\begin{equation}\label{2.142}
  \lim_{z\rightarrow 1-} \frac{V^*(\lambda-\lambda z)-1}{z-1}=\lambda
  E[V],
\end{equation}
and
\begin{eqnarray}\label{2.143}
  \lim_{z\rightarrow 1-}\frac{z-S^*(\lambda-\lambda z)}{(1-z)^2 \ln(\lambda-\lambda z)}
  &=& \lim_{z\rightarrow
  1-}\frac{1-{S^*}'(\lambda-\lambda
  z)}{(1-z)[2 \ln(\lambda-\lambda z)-(1-z)]}\nonumber \\
   &=& \lim_{z\rightarrow
  1-}\frac{-{S^*}''(\lambda-\lambda
  z)}{2 \ln(\lambda-\lambda z)+2(1-z)+2}\nonumber \\
  &=& \frac{c \lambda^2}{2}.
\end{eqnarray}
It follows from (\ref{2.142}) and (\ref{2.143})  that
\begin{eqnarray*}
 \lim_{z\rightarrow 1-}\Pi(z) (1-z)\ln \left (\frac{1}{\lambda - \lambda z}\right )
&=& \frac{\pi_0}{(1-\nu_0)} \lim_{z\rightarrow
 1-} \frac{[V^*(\lambda-\lambda z)-1]S^*(\lambda-\lambda z)}{z-S^*(\lambda-\lambda z)} (1-z)\ln \left (\frac{1}{\lambda - \lambda z}\right ) \nonumber\\
  &=&  \frac{2 \pi_0 E[V]}{c(1-\nu_0)
  \lambda}.
\end{eqnarray*}
Hence, we have
\[
  \Pi(z) \sim \frac{1}{(1-z)\ln (\frac{1}{\lambda - \lambda z})}\sim
  \frac{2 \pi_0 E[V]}{c(1-\nu_0)
  \lambda}
  \frac{1}{(1-z)\ln (\frac{1}{1 - z})}, \ \ \mbox{as}\ \ z\rightarrow 1-.
\]
Using another Tauberian theorem (Theorem~\ref{the:4.3})  leads to
\[
  S_{\pi}(n) \sim \frac{2 \pi_0  E[V]}{c(1-\nu_0)
  \lambda}\frac{n}{\ln n},
\]
from which and Lemma~\ref{lem:2.1}, (\ref{2.14}) follows
immediately. \BOX

\begin{remark} (a) In the assumption of  $s(x)\sim c x^{-\theta}$ in case (ii),
it is reasonable to restrict $\theta$ in the domain $(2, 3]$.
Otherwise, if $\theta =2$, then $E[S]=\infty$, and if $\theta>3$,
then $E[S^2]<\infty$, both contradict other assumptions. (b) The
condition that $E[S^2]<\infty$ is often imposed in studies of this
type of problems, for example in \cite{Bai1992}, but as far as we
know, Case (ii) is a new one. (c) In the case $\rho=1$, the
asymptotic behavior of the loss probability is independent of the
vacation time $V$.
\end{remark}

\subsection{The case when $\rho>1$}

This is the simplest case, in which we no longer need to impose
any further property assumptions on the service time or on the
vacation time for the following asymptotic property for the loss
probability.

\begin{theorem} \label{the:2.4} Assume that  $\rho>1$. Then we have
\begin{equation}\label{2.15}
  \lim_{N\rightarrow \infty} z_2^{-N}\left[P_{\mbox{loss}}(N) - \left (1-\frac{1}{\rho} \right )\right]
  =\frac{(1-z_2)[1-{S^*}'(\lambda-\lambda z_2)]E[V]}{ z_2 \rho E[S] [1-V^*(\lambda-\lambda z_2)]},
\end{equation}
 where $z_2\in (0,1)$ is the unique solution
 to the equation $z=S^*(\lambda-\lambda z)$.
\end{theorem}

\proof From Lemma 2.2, we know that there exists a unique solution
$z_2\in (0,1)$ to the equation $z=S^*(\lambda-\lambda z)$.
Obviously, ${S^*}'(\lambda-\lambda z_2)<1$. Since
\begin{eqnarray*}
\lim_{z\rightarrow z_2-} \Pi(z)(z_2-z) = \frac{z_2\pi_0
[1-V^*(\lambda-\lambda z_2)] ]}{(1-\nu_0) [1-{S^*}'(\lambda-\lambda
z_2)]},
\end{eqnarray*}
 $z_2$ is a simple pole for $\Pi(z)$. Similar to the proof of Theorem~\ref{the:2.2}, we know that
$z=z_2$ is the unique singular point of $\Pi(z)$ on the circle
$|z|=z_2$. By Theorem~5.2.1 in \cite{wilf94}, a standard result on
asymptotics of complex functions, we have
\[
\pi_n\sim \frac{\pi_0 [1-V^*(\lambda-\lambda z_2)]}{
(1-\nu_0)[1-{S^*}'(\lambda-\lambda z_2)]}z_2^{-n},
\]
from which and the Stolz-Ces\'{a}ro theorem it follows that
\begin{equation*}
S_{\pi}(n)\sim \frac{z_2 \pi_0 [1-V^*(\lambda-\lambda z_2)]
}{(1-z_2)(1-\nu_0) [1-{S^*}'(\lambda-\lambda z_2)]}z_2^{-n}.
\end{equation*}
The assertion (\ref{2.15}) follows mow from Lemma~\ref{lem:2.1}
immediately.\BOX

\section{Decay for standard M/G/1/N and GI/M/1/N queues}

The M/G/1/N queue with vacations becomes the (standard) M/G/1/N
queue when the vacation time $V=0$. The embedded transition matrix
of the  M/G/1/N queue is obtained by replacing  $b_k$ by $a_k$,
$k=0,\cdots, N-2$, in the first row of the embedded transition
matrix for the M/G/1/N queue with vacations. Denote by
$\hat{P}_{\mbox{loss}}(K)$ the loss probability for the M/G/1/N
queue. It is known that
\[
  \hat{P}_{\mbox{loss}}(N) =
  1-\frac{1}{\hat{\pi}_0(N)+\rho},
\]
where $\hat{\pi}(N)$ is the invariant vector of the embedded
Markov chain for the $M/G/1/N$ queue. The following result is
obtained by the same asymptotic analysis for $P_{\mbox{loss}}(N)$
in Section 2 with $V= 0$.

\begin{theorem}\label{the:3.1} For the M/G/1/N queue, we have the
following asymptotic results on the loss probability
$\hat{P}_{\mbox{loss}}(N)$.
\begin{description}
\item[(i)] Let $\rho<1.$

(1) If $R_{S^*}=1$
  and $\overline{S}(x)\sim x^{-\alpha}
L(x)$ for some $\alpha>1$, then
\begin{equation*}
  \lim_{N\rightarrow \infty} \frac{N^{\alpha-1}}{L(\frac{N}{\lambda})} \hat{P}_{\mbox{loss}}(N)
  =\frac{\lambda^\alpha}{\alpha-1}.
\end{equation*}

(2) If $R_{S^*}>1$ and  $R_{S^*}< S^*(\lambda-\lambda R_{S^*})$,
then there exists only one solution  $\sigma_1$ to the equation
$z=S^*(\lambda-\lambda z)$, $1<z< R_{S^*}$ such that
\begin{equation*}
  \lim_{N\rightarrow \infty} \sigma_1^{N} \hat{P}_{\mbox{loss}}(N)
  =\frac{ \sigma_1 (1-\rho)^2}{[{S^*}'(\lambda-\lambda \sigma_1)-1]}.
\end{equation*}

\item[(ii)] Let  $\rho=1$.

(1) If $E[S^2]<\infty$, then
\begin{equation*}
  \lim_{N\rightarrow \infty} N \hat{P}_{\mbox{loss}}(N)=\frac{\lambda^2 E[S^2]}{2}.
\end{equation*}

(2) If $E[S^2]=\infty$ and $s(x)\sim c x^{-\theta}$, then for
$2<\theta<3$,
\begin{equation*}
 \lim_{N\rightarrow \infty} N^{\theta-2} \hat{P}_{\mbox{loss}}(N)=\frac{c \lambda^{\theta-1}
  \Gamma(\theta-1)\Gamma(4-\theta)}{(1-\theta)(2-\theta)(3-\theta)},
\end{equation*}
and for $\theta=3$,
\begin{equation}\label{2.141b}
 \lim_{N\rightarrow \infty} \frac{N}{\ln N} \hat{P}_{\mbox{loss}}(N)=\frac{c
 \lambda^2 }{2}.
\end{equation}

\item[(iii)] Let  $\rho>1$.
 then
\begin{equation*}
  \lim_{N\rightarrow \infty} \sigma_2^{-N}\left[\hat{P}_{\mbox{loss}}(N) - \left (1-\frac{1}{\rho} \right )\right]
  =\frac{1-{S^*}'(\lambda-\lambda \sigma_2)}{ \sigma_2 \rho^2 },
\end{equation*}
where $\sigma_2\in (0,1)$ is the unique solution
 to the equation $z=S^*(\lambda-\lambda z)$.
\end{description}
\end{theorem}

\begin{remark} (a) The above results for $\rho>1$ and for $\rho=1$
with the assumption that $E[S^2]<\infty$ coincide with those in
(\cite{Bai1992}). (b) When $\rho <1$, we show that the asymptotic
rate can be exactly geometric or polynomial (heavy-tailed). The
later case was neglected in \cite{Bai1992}. (c) The asymptotic
result for the case that $\rho=1$ with $E[S^2]=\infty$ is considered
new.
\end{remark}

We now consider the (standard) $GI/M/1/N$ queue.  Denote by $A$
the generic interarrival time random variable, and by $A(x)$ and
$a(x)$ its distribution function and probability density function,
respectively. Let $\mu$ be the service rate and let $N$ be the
maximal capacity of the system, or the number of customers in the
system including the one in service. The traffic intensity of the
$GI/M/1/N$ queues is $\tilde{\rho}=\frac{1}{\mu E[A]}$. By
$\tilde{P}_{\mbox{loss}}(N)$ we denote the loss probability of the
GI/M/1/N queue. Using the dual concept (see, e.g.
\cite{{Bai1992},{Miy90}}), we have
\begin{equation}\label{3.1}
{\tilde P}_{\mbox{loss}}(N)=\tilde{\pi}_N (N)=\hat{\pi}_0(N+1),
\end{equation}
where $\tilde{\pi}(N)$ is the invariant distribution of the
embedded $GI/M/1/N$ queue, and $\hat{\pi}(N+1)$ is the invariant
distribution of the $M/G/1/N+1$ queue in which the arrival process
 is Poissonian with parameter $\mu$, the distribution of the service
time $A(x)$ and the traffic intensity $\rho=\tilde{\rho}^{-1}$. From
(\ref{2.4}), (\ref{3.1}) and the proof of
Theorems~\ref{the:2.1}--\ref{the:2.3}, we obtain the following
assertion.

\begin{theorem}\label{the:3.2} Consider the GI/M/1/N queue.
\begin{description}
\item[(i)] Let $\tilde{\rho}<1$. We then have
\begin{equation*}
  \lim_{N\rightarrow \infty} \eta_1^{-N} \tilde{P}_{\mbox{loss}}(N)
  = 1-{A^*}'(\mu -\mu \eta_1),
\end{equation*}
where $\eta_1 \in (0,1)$ is the unique solution
 to the equation $z=A^*(\mu-\mu z)$.

\item[(ii)] Let  $\tilde{\rho}=1$.

(1) If $E[A^2]<\infty$, then
\begin{equation*}
  \lim_{N\rightarrow \infty} N \tilde{P}_{\mbox{loss}}(N)=\frac{\mu^2 E[A^2]}{2}.
\end{equation*}

(2) If $E[A^2]=\infty$ and $a(x)\sim c x^{-\theta}$,  then for
$2<\theta< 3$,
\begin{equation*}
 \lim_{N\rightarrow \infty} N^{\theta-2} \tilde{P}_{\mbox{loss}}(N)=\frac{c \mu^{\theta-1}
  \Gamma(\theta-1)\Gamma(4-\theta)}{(1-\theta)(2-\theta)(3-\theta)},
\end{equation*}
and for $\theta =3$,
\begin{equation}\label{2.141c}
 \lim_{N\rightarrow \infty} \frac{N}{\ln N} \tilde{P}_{\mbox{loss}}(N)=\frac{c
 \mu^2 }{2}.
\end{equation}

\item[(iii)] Let  $\tilde{\rho}>1$ and let $R_{A^*}$ be the
leftmost singular point of the function $A(\mu-\mu z)$.

 (1) If $R_{A^*}=1$
  and $\overline{A}(x)\sim x^{-\alpha}
L(x)$ for some $\alpha>1$ as $x\rightarrow \infty$, then
\begin{equation*}
  \lim_{N\rightarrow \infty} \frac{N^{\alpha-1}}{L(\frac{N}{\lambda})}
  \left[\tilde{P}_{\mbox{loss}}(N)-\frac{1}{\tilde{\rho}}\right]
  =\frac{\lambda^\alpha}{\alpha-1}.
\end{equation*}
(2) If $R_{A^*}>1$ and $R_{A^*}< S^*(\lambda-\lambda R_{A^*})$, then
there exists only one solution $\eta_2$ to the equation
$z=A^*(\mu-\mu z)$, $z\in (1,R_{A^*})$ such that
\begin{equation*}
  \lim_{N\rightarrow \infty} \eta_2^{N} \left[\tilde{P}_{\mbox{loss}}(N)-\frac{1}{\tilde{\rho}}\right]
  =\frac{(1-\rho)^2}{A^*(\mu-\mu \eta_2) -1}.
\end{equation*}
\end{description}
\end{theorem}

\begin{remark} (a) The results for $\tilde{\rho}<1$ and for $\tilde{\rho}=1$
with $E[A^2]<\infty$ coincide with those in
\cite{Bai1992,{Kim2000}}. (b) When $\tilde{\rho}>1$, we showed that
the asymptotic rate can be either geometric or polynomial. The
latter case was not considered in (\cite{Bai1992}). (c) The
asymptotic result on the case that $\tilde{\rho}=1$ with
$E[A^2]=\infty$ is considered new.
\end{remark}

\section{Conclusions}

The main contribution of this paper is a detailed characterization of the
asymptotic rates of the loss probability for the M/G/1/N-(V,E)-queue. New asymptotic behavior was found even for simpler  M/G/1/N queue and the GI/M/1/N queue in terms of introduction of new asymptotic analysis methods.

Although we only focused on asymptotics of the loss probability, the arguments in this paper can be extended to derive asymptotic results of the other performance measures, such as the queue length distribution.  Let $L(t)$ be the queue length process of the M/G/1/$\infty$-(V, E)-queue. Suppose that $L(t)$ is stable, i.e. $\rho<1$. It is well known that $L(t)$ has the invariant distribution $\pi$ given by (2.3), i.e. $L(t)$ has the same invariant distribution as that of its embedded queue. Define
\[
  \|\pi^*(N) - \pi \| = \sum_{j=0}^N|\pi_j^*(N)-\pi_j| + \sum_{j=N+1}^\infty \pi_j,
\]
where $\pi^*_N(N)$ is given by (2.2) and $\pi^*_j(N)$ is given by (see formula (13) in \cite{Fr97})
\[
  \pi^*_j(N) = \frac{\pi_j(N)(1-\nu_0)\lambda^{-1}}{E[V]\pi_0(N)+E[S](1-\nu_0)},\ \ 0\leq j\leq  N-1.
\]
Since $\pi_j(N)=\frac{\pi_j}{S_{\pi}(N)}$, $0\leq j\leq  N-1$, and $\pi_0=\frac{(1-\rho)(1-\nu_0)}{\lambda E[V]}$, it is not hard to find that the asymptotic rate of  $\|\pi^*(N) - \pi \|$ as $N\rightarrow \infty$ is completely determined by the decay rate of $\overline{S}_{\pi}(N)$ as $N\rightarrow \infty$. Hence the asymptotic rate of $\|\pi^*(N) - \pi \|$ can be derived by using similar arguments to that in Section 3.1.

\begin{center}
\bigskip \noindent{\bf\large Acknowledgements}
\end{center}

Both authors are deeply grateful to the anonymous referee for very constructive comments, which have allowed the authors to
improve the presentation of this paper significantly. This work was supported in part by an NSERC discovery grant of Canada, the Fundamental Research Funds for the Central Universities (grant number 2010QYZD001), and the National Natural Science Foundation of China (grant numbers  10901164, 11071258).

\appendix

\section{Appendix}


 The following Definition A.1  and Theorem A.1 are taken from
Definition VI.1 and Corollary VI.1 in \cite{Fla09}, respectively, with a slight
difference. Theorem A.1 is called Tauberian-like theorem in \cite{Li2011}.

\begin{definition} Let $r\geq 1$.
For given constants $\phi, \hat{r}$ with $\hat{r}>r$ and
$0<\beta<\frac{\pi}{2}$, define $\Delta(r, \beta, \hat{r})$ to be
the following open domain
\[
 \Delta(r, \beta, \hat{r})=\{z:|z|<\hat{r}, z\neq r, |Arg(z-r)|>\beta \}.
\]
A domain is called a $\Delta(r)$-domain if it is a $\Delta(r, \beta,
\hat{r})$ for some $\hat{r}$ and $\beta$. A function is said to be
$\Delta(r)$-analytic if it is analytic in some $\Delta(r)$-domain.
\end{definition}

\begin{theorem} \label{the:4.1}
 Assume that $f(z)$ is $\Delta(1)$-analytic and
\[
 f(z)\sim \frac{c}{(1-z)^\theta},\ \ \ \mbox{as} \ \ z\rightarrow 1,\ \ z\in \Delta(1)
 \]
 for some constant $c$ and $\theta$.  Let $f_n$ be  the $n$th Taylor coefficient of
$f(z)$. (i) If $\theta \notin \{0,-1,-2,\cdots\}$, then
\[
 f_n \sim \frac{c }{\Gamma(\theta)}n^{\theta-1};
\]
(ii) If  $\theta \in \{0,-1,-2,\cdots\}$, then
\[
 f_n = o(n^{\theta-1}).
\]
\end{theorem}

The following theorems, referred to as Tauberian theorems, are taken
from Section~5 of Chapter~13 in \cite{Feller71}.

\begin{theorem} \label{the:4.2}
Let $U$ be a measure on $[0, \infty]$ with Laplace transform
$\omega(\tau)=\int_0^\infty e^{-\lambda x}U\{dx\}$.
 If $L$ is
a slowly varying function and $0\leq \rho<\infty$, then the
properties
\[
\omega(\tau) \sim \tau^{-\rho} L\left (\frac{1}{\tau}\right ),\ \
\tau \rightarrow 0,
\]
and
\[
    U(t)\sim \frac{1}{\Gamma(\rho+1)} t^{\rho} L(t),\ \ t\rightarrow
    \infty,
\]
are equivalent.
\end{theorem}

\begin{theorem}\label{the:4.3}
Let $q_n\geq 0$ and assume that
\[Q(s)=\sum_{n=0}^\infty q_n s^n\]
converges for $0\leq s<1$. If $L$ is a slowly varying function and
$0\leq \rho<\infty$, then the properties
\begin{equation}\label{7.1}
    Q(s)\sim \frac{1}{(1-s)^{\rho}}L\left (\frac{1}{1-s}\right ),\ \ s\rightarrow 1-,
\end{equation}
and
\[
    q_0+q_1+\cdots+q_{n-1}\sim \frac{1}{\Gamma(\rho+1)}n^{\rho}L(n),\ \ n\rightarrow \infty,
\]
are equivalent.

\end{theorem}

\end{document}